\documentclass{article}

\def\ra{\rightarrow}

\def\ss{\subseteq}

\def\Re{\hbox{\rm Re}\,}

\def\K{{\cal K}}
\def\C{{\cal C}}
\def\M{{\cal M}}

\def\avgint{\frac{1}{2\pi} \int_0^{2\pi}}
\def\dbar{\overline{\partial}}

 \def\HollowBox #1#2{{\dimen0=#1 \advance\dimen0 by -#2       
       \dimen1=#1 \advance\dimen1 by #2                       
        \vrule height #1 depth #2 width #2                    
        \vrule height 0pt depth #2 width #1                   
        \llap{\vrule height #1 depth -\dimen0 width \dimen1}%
       \hskip -#2                                             
       \vrule height #1 depth #2 width #2}}                   
 \def\BoxOpTwo{\mathord{\HollowBox{6pt}{.4pt}}\;}             

\def\endpf{\hfill $\BoxOpTwo$}

\font\teneufm=eufm10
\font\seveneufm=eufm7
\font\fiveeufm=eufm5
\newfam\eufmfam
\textfont\eufmfam=\teneufm
\scriptfont\eufmfam=\seveneufm
\scriptscriptfont\eufmfam=\fiveeufm

\newfam\msbfam
\font\tenmsb=msbm10  \textfont\msbfam=\tenmsb
\font\sevenmsb=msbm7  \scriptfont\msbfam=\sevenmsb
\font\fivemsb=msbm5    \scriptscriptfont\msbfam=\fivemsb
\def\Bbb{\fam\msbfam \tenmsb}

\def\O{\Omega}

\def\CC{{\Bbb C}}

\newtheorem{theorem}{Theorem}
\newtheorem{corollary}[theorem]{Corollary}
\newtheorem{proposition}[theorem]{Proposition}
\newtheorem{lemma}[theorem]{Lemma}
\newtheorem{remark}[theorem]{Remark}

\newtheorem{definition}{Definition}
\newtheorem{example}[definition]{EXAMPLE}

\makeindex

\begin{document}

\begin{center}
\huge \bf
The Kobayashi Metric, \\
Extremal Discs, \\
and Biholomorphic Mappings\footnote{{\bf Key Words:}  Kobayashi metric,
extremal disc, Carath\'{e}odory metric, pseudoconvexity.}\footnote{{\bf MR Classfication Numbers:}
32E05, 32H02, 32T27, 32T05.}
\end{center}
\vspace*{.12in}

\begin{center}
\large Steven G. Krantz\footnote{Author supported in part
by the National Science Foundation and the Dean of the Graduate
School at Washington University.}
\end{center}
\vspace*{.15in}

\begin{center}
\today
\end{center}
\vspace*{.2in}

\begin{quotation}
{\bf Abstract:} \sl
We study extremal discs for the Kobayashi metric.  Inspired by work of
Lempert on strongly convex domains, we present results on strongly pseudoconvex
domains.

We also consider a useful biholomorphic invariant, inspired by the Kobayashi
(and Carath\'{e}odory) metric, and prove several new results about biholomorphic
equivalence of domains.  Some useful results about automorphism groups of
complex domains are also established.
\end{quotation}
\vspace*{.25in}

\setcounter{section}{-1}

\section{Introduction}

Throughout this paper, a {\it domain} in $\CC^n$ is a connected, open set.  Usually
our domains will be bounded.   It is frequently convenient to think
of a domain $\Omega$ (with smooth boundary) as given by
$$
\Omega = \{z \in \Omega: \rho(z) < 0\} \, ,
$$
where $\rho$ is a $C^k$ function and $\nabla \rho \ne 0$ on $\partial \Omega$.
We say in this circumstance that $\rho$ is a $C^k$ {\it defining function} for $\Omega$.
It follows from the implicit function theorem that $\partial \Omega$ is a $C^k$ manifold
in a natural sense.   See [KRA1] for more on these matters.

Throughout the paper $D$ denotes the unit disc in the complex plane $\CC$ and $B$ denotes
the unit ball in complex space $\CC^n$.
If $\Omega_1, \Omega_2$ are domains in complex space then we let
$\Omega_1(\Omega_2)$ denote the holomorphic mappings from $\Omega_2$
to $\Omega_1$.   In case $\Omega_2$ is either $D$ or $B$ and $z \in \Omega$ then we sometimes
let $\Omega^z(D)$ (resp.\ $\Omega^z(B)$) denote the elements $\varphi \in \Omega(D)$ (resp.\ $\varphi \in \Omega(B)$)
such that $\varphi(0) = z$.

The infinitesimal {\it Kobayashi metric} on $\Omega$ is defined as follows.  Let $z \in \Omega$
and $\xi \in \CC^n$.  Then
\begin{eqnarray*}
F_\Omega^K(z, \xi) & = & \inf\{\alpha: \alpha > 0\ \mbox{\rm and} \ \exists 
f \in \Omega(D)\ \mbox{\rm with} \ f(0) = z,  f'(0) = \xi/\alpha \} \\
    & = & \inf \left\{\frac{|\xi|}{|f'(0)|}: f \in \Omega^z (D) \right \} \, .
\end{eqnarray*}
The infinitesimal {\it Carath\'{e}odory metric} is given by
$$
F_C^\Omega(z,\xi) \equiv \sup_{f \in D(\O) 
                       \atop
                       f(z) = 0}
                                  |f'(z) \xi| \, .
$$
In these definitions, $| \ \ |$ denotes Euclidean length.
The definitions of both these metrics are motivated by the proof of the Riemann mapping theorem, and by
the classical Schwarz lemma.  Details may be found in [KRA1] and [KRA2].

Companion notions are the Kobayashi and Carath\'{e}odory volume elements.  We define
these as follows (see also [EIS]).   If $\Omega$ is a fixed domain and $z \in \Omega$
then set
$$
{\cal C}_\Omega(z) = {\cal C}(z) = \sup\{|\hbox{det} \, \varphi'(z)| : \varphi: \Omega \rightarrow B, \varphi (z) = 0 \} 
$$
and
$$
{\cal K}_\Omega(z) = {\cal K}(z) = \inf \left \{ \frac{1}{|\hbox{det} \, \psi'(z)|} : \psi: B \rightarrow \Omega, \psi(0) = z \right \} \, .
$$

If $\varphi$ is a candidate mapping for ${\cal C}$ and $\psi$ is a candidate mapping for ${\cal K}$, then an examination
of $\varphi \circ \psi$ using the Schwarz lemma (see [RUD]) shows that ${\cal C}(z) \leq {\cal K}(z)$ for any $z \in \Omega$.
We set 
$$
{\cal M}(z) = \frac{{\cal K}(z)}{{\cal C}(z)} \,  .
$$
We call ${\cal M}$ the {\it quotient invariant}.  Of course ${\cal M}(z) \geq 1$ for all $z \in \Omega$.  The following remarkable lemma
of Bun Wong (see [WON]) is useful in the study of automorphism groups:

\begin{lemma} \sl
Let $\Omega \ss \CC^n$ be a bounded domain.  If there is a point $z \in \Omega$ so that ${\cal M}(z) = 1$ then
$\Omega$ is biholomorphic to the unit ball $B$ in $\CC^n$.
\end{lemma}

We shall not prove this result here, but refer the reader instead to [KRA1].  
It is worth stating the fundamental result of Bun Wong and Rosay (again
see [KRA1] for the details) that is proved using Lemma 1.

\begin{theorem} \sl
Let $\Omega \ss \CC^n$ be a bounded domain and $P \in \partial \Omega$ a point of
strong pseudoconvexity.  Fix a point $X \in \Omega$ and suppose that there
are biholomorphic mappings $\varphi_j: \Omega \ra \Omega$ ({\it automorphisms} of $\Omega$) so
that $\varphi_j(X) \ra P$ as $j \ra \infty$.  Then $\Omega$ is biholomorphic to the 
unit ball $B$ in $\CC^n$.
\end{theorem}

This theorem has been quite influential in the development of the theory of automorphism
groups of smoothly bounded domains.  See, for example, [GRK1], [ISK], and [GKK].
It is common to call the point $P$ in the theorem a {\it boundary orbit accumulation
point for the automorphism group action} (or ``orbit accumulation point'' for short).

\section{The Quotient Invariant}

Here we discuss in detail the invariant of Bun Wong described in Section 0.
It has far-reaching implications beyond the basic application in the proof
of the Bun Wong/Rosay theorem.

\begin{proposition} \sl  Let $\Omega \ss \CC^n$ be a bounded domain.  If there
is a point $P \in \Omega$ such that ${\cal M}(P) = 1$ then
${\cal M}(z) = 1$ for all $z \in \Omega$.  Obversely, if there
is a point $P \in \Omega$ with ${\cal M}(P) > 1$ then ${\cal M}(z) > 1$
for all $z \in \Omega$.
\end{proposition}
{\bf Proof:}  If ${\cal M}(P) = 1$ for some $P$ then Bun Wong's original lemma
 (Lemma 1 above) shows that $\Omega \cong B$.  Of course ${\cal M}$ is a biholomorphic
 invariant.  And $B$ has transitive automorphism group.  It follows therefore
 that $\Omega$ has invariant ${\cal M}$ with value 1 at every point.

Obversely, if ${\cal M}(P) \ne 1$ at some point then, by contrapositive reasoning
in the last paragraph, it cannot be that ${\cal M}$ equals 1 at any point.

That completes the proof of the proposition.
\endpf
\smallskip \\

\begin{proposition} \sl
Let $\Omega \ss \CC^n$ be a bounded domain.  Let $P \in \partial \Omega$ and suppose
that $\partial \Omega$ is $C^2$ and strongly  pseudoconvex near $P$.
Then
$$
\lim_{z \ra P} {\cal M}(z) = 1 \, .
$$
\end{proposition}
{\bf Proof:}  This follows from the asymptotics of I. Graham for the Carath\'{e}odory and
Kobayashi metrics on such a domain.  The main point is that $\partial \Omega$ is approximately
a ball near $P$, so the asymptotic behavior of $F_K^\Omega$, $F_C^\Omega$, ${\cal K}_\Omega$,
and ${\cal C}_\Omega$ is the same as that on the domain the ball $B$.  
\endpf
\smallskip \\

\begin{proposition} \sl
Let 
$$
{\cal E} = \{(z_1, z_2, \dots, z_n) \in \CC^n: |z_1|^{2m_1} + |z_2|^{2m_2} + \cdots + |z_n|^{2m_n} < 1\}
$$
be a domain in $\CC^n$, with $m_1, m_2, \dots, m_n$ positive integers.   Often $E$ is called
an {\it egg} or an {\it ellipsoid}.   If some $m_j > 1$ then
$E$ is {\it not} biholomorphic to the ball.
\end{proposition}
{\bf Proof:}   This result was first proved by S. Webster [WEB] using techniques of differential
geometry.  Later, S. Bell [BEL] gave a very natural proof by showing that any biholomorphism
of the ball to $E$ must extend smoothly to the boundary, and then noting that the Levi
form is a biholomorphic invariant.  Here we give a proof that uses $\M$.

For simplicity we shall take $n = 2$, $m_1 = 1$, and $m_2 > 1$.  Seeking a contradiction, we let $\varphi: B \ra E$ be
a holomorphic mapping that takes 0 to 0.  Thus $\varphi$ is a candidate mapping for the calculation
of $\K_\Omega$.  Now set
$$
\widetilde{\varphi}(z_1, z_2) = \frac{1}{4\pi^2} \int_0^{2\pi} \!\!\!\int_0^{2\pi} f(z_1 e^{i\theta_1}, z_2 e^{i\theta_2}) e^{-i\theta_1} e^{-i\theta_2} \, d\theta_1 d\theta_2 \, .
$$
Then one may calculate that {\bf (i)}  $\widetilde{\varphi}$ still maps $B$ into $E$ and {\bf (ii)}  the first (holomorphic) derivatives
of $\widetilde{\varphi}$ at 0 are the same as the first (holomorphic) derivatives of $\varphi$ at 0.
Also $\widetilde{\varphi}$ is linear (since the higher-order terms all average to 0).

As a result of the last paragraph, we may calculate $\K$ at 0 for $E$ using only linear maps.  A similar
argument applies to maps $\psi: E \ra B$.  Of course it is obvious that there is no linear
equivalence of $B$ and $E$ (the boundaries of the two domains have different curvatures, for instance).  In particular, $\M(0) > 1$.  It follows that $\M(P) > 1$ at all points
$P$ of $E$.   Thus $E$ and $B$ are biholomorphically inequivalent.
\endpf
\smallskip \\

\begin{proposition} \sl
Let $\Omega \ss \CC^n$ be a bounded domain with $C^2$ boundary.  If $P \in \partial \Omega$ is a point
of strong pseudoconcavity, let $\nu$ be the unit outward normal vector at $P$.  Set $P_\epsilon = P - \epsilon \nu$.
Then $\M_\Omega(P_\epsilon) \approx C \cdot \epsilon^{-3/4}$.
\end{proposition}
{\bf Proof:}  It is a result of [KRA3] that the Kobayashi metric $F_K^\Omega(P_\epsilon, \nu)$ is of
size $C \cdot \epsilon^{-3/4}$.  It is also clear that the Kobayashi metric at $P_\epsilon$ in
complex tangential directions is of size $C$, where $C > 0$ is some universal positive constant.
Hence ${\cal K} \sim C \cdot \epsilon^{-3/4}$.  On the other hand, the Hartogs extension phenomenon
gives easily that ${\cal C}(P) \sim C$.  It follows then that $\M \approx C \cdot \epsilon^{-3/4}$.
\endpf
\smallskip \\

\begin{corollary} \sl
Let $\Omega \ss \CC^n$ be a bounded domain with $C^2$ boundary.  If $P \in \partial \Omega$ is a point
of strong pseudoconcavity, then $P$ cannot be a boundary orbit accumulation point.
\end{corollary}
{\bf Proof:}   Seeking a contradiction, we suppose that $P$ is a boundary orbit accumulation point.
So there is a point $X \in \Omega$ and there are automorphisms $\varphi_j$ of $\Omega$ so
that $\varphi_j(X) \ra P$.  But of course $\M(X)$ is some positive constant $C$ that exceeds 1.
And the invariant $\M(z)$ blows up like $\hbox{dist}(z, \partial \Omega)^{-3/4}$ as $z \ra P$.
This is impossible.
\endpf
\smallskip \\

\begin{remark} \rm
It is a result of [GRK2] that if $\Omega$ is {\it any} domain and $P \in \partial \Omega$ a point
of non-pseudoconvexity (even in the weak sense of Hartogs) then $P$ cannot be a boundary orbit
accumulation point.  The last Corollary captures a special case of this result
using the idea of the quotient invariant.
\end{remark}

\begin{proposition} \sl
Let $\Omega \ss \CC^2$ be a smoothly bounded domain that is of finite type (in the sense of Kohn/D'Angelo/Catlin---see [KRA1])
at every boundary point.  Let $P \in \partial \Omega$.  Then 
$$
0 < C_1 \leq \liminf_{z \ra P} \M(z) \leq \limsup_{z \ra P} \M(z) \leq C_2
$$
for some universal, positive constants $C_1$, $C_2$.
\end{proposition}
{\bf Proof:}   This follows from the estimates in [CAT].
\endpf 
\smallskip \\

\begin{proposition} \sl
Let $\Omega \ss \CC^2$ be a smoothly bounded, convex domain of finite type.   Let $P \in \partial \Omega$.
Then
$$
0 < C_1 \leq \liminf_{z \ra P} \M(z) \leq \limsup_{z \ra P} \M(z) \leq C_2
$$
for some positive constants $C_1$, $C_2$.
\end{proposition}
{\bf Proof:} Fix a point $z \in \Omega$ near $P$ and $\xi$ a
tangent direction at $z$. Certainly any mapping $\varphi: D
\ra \Omega$, $\varphi(0) = z$ with $\varphi'(0) = \lambda \xi$
for some $\lambda > 0$ is a candidate for the Kobayashi metric
at $z$ in the direction $\xi$, and the reciprocal of its
derivative gives an upper bound for the Kobayashi metric at
that point in that direction. In particular, we may take
$\varphi$ to be the obvious linear embedding of the disc $D$
into $\Omega$ pointing in the direction $\xi$ (with image having
diameter $\delta$, the distance from $z$ to $\partial \Omega$ in the direction $\xi$) and with
$\varphi(0) = z$.

Thanks to work of McNeal [MCN], we know that the type of a convex point of finite type
can be measured with the order of contact by complex lines.  If, after a rotation
and translation, we take $P$ to be the point $(1,0)$ and $\langle 1, 0\rangle$ the real
normal direction, then the complex line of greatest contact will of course be
$\zeta \mapsto (1, \zeta)$.  Let that order of contact be $2m$ for some positive integer $m$.
Then it is clear, after shrinking $\Omega$ if necessary, that an ellipsiod of the form
$$
E = \{(z_1, z_2) \in \CC^2: |z_1|^2 + K |z_2|^{2m} < 1\}  
$$
will osculate $\partial \Omega$ at $(1,0)$ and will contain $\Omega$.  So, in particular
$F_C^\Omega(z, \xi) \geq F_C^E(z, \xi)$ for any $z \in \Omega$ and $\xi$ any tangent
vector.

We calculate that, for $z = (\alpha, 0) \in E$, the mappings
$$
(\zeta_1, \zeta_2) \longmapsto \frac{\zeta_1 - \alpha}{1 - \overline{\alpha}z_1}
$$
and
$$
(\zeta_1, \zeta_2) \longmapsto \frac{\sqrt[2m]{1 - |\alpha|^2} z_2}{1 - \overline{\alpha}z_1}
$$
are candidate maps for the Carath\'{e}odory metric at the point $z$.  The first one gives
a favorable lower bound for the Carath\'{e}odory metric in the normal direction $\langle 1, 0\rangle$
at $z$ and the second gives a favorable lower bound for the Carath\'{e}odory metric in the tangential
direction $\langle 0, 1\rangle$ at $z$.   Of course these are also lower bounds for the Carath\'{e}odory metric
on $\Omega$.

It is easy to see that the given upper bounds for the Kobayashi metric and the given lower
bounds for the Carath\'{e}odory metric are comparable.  Since $F_C^\Omega \leq F_K^\Omega$ always
(see [KRA1]), it follows that $\M \approx C$ (a constant) on a smoothly bounded, convex domain of finite type in $\CC^2$.
\endpf
\smallskip \\

\begin{remark} \rm
The elementary comparison of the domains $\Omega$ and $E$ that we exploited in the last
proof will not work in higher dimensions.  The matter in that context is more
subtle.
\end{remark}

\section{More on the Quotient Invariant}

It is natural to wonder about the role of the ball $B$ in the definition of the quotient invariant
$\M$.  We define $\K$ in terms of mappings from the ball $B$ to the given domain $\Omega$ and
we define $\C$ in terms of mappings from the given domain $\Omega$ to the ball $B$.  What
if the ball $B$ were to be replaced by some other ``model domain''?

Let ${\cal B}$ be some fixed, bounded domain in $\CC^n$.  Fix a point $P_0 \in {\cal B}$.
Let $\Omega$ be some other bounded domain, and let $z \in \Omega$.  Define new invariants
$$
\widehat{\cal C}_\Omega(z) = \widehat{\cal C}(z) = \sup\{|\hbox{det} \, \varphi'(z)| : \varphi: \Omega \rightarrow {\cal B}, \varphi (z) = P_0 \} 
$$
and
$$
\widehat{\cal K}_\Omega(z) = \widehat{\cal K}(z) = \inf \left \{ \frac{1}{|\hbox{det} \, \psi'(z)|} : \psi: {\cal B} \rightarrow \Omega, \psi(P_0) = z \right \} 
$$
and a new quotient invariant
$$
\widehat{\M}_\Omega(P) = \widehat{\M}(P) = \frac{\widehat{\K}_\Omega(P)}{\widehat{\C}_\Omega(P)} \, .
$$

Now we have

\begin{proposition} \sl
Let $\Omega$ be any given bounded domain in $\CC^n$.  Suppose that there is a point
$P \in \Omega$ such that $\widehat{\M}_\Omega(P) = 1$.  Then
$\Omega$ is biholomorphic to the model domain ${\cal B}$.
\end{proposition}
{\bf Proof:}  The argument is just the same as in the classical case of ${\cal B} = B$, the
unit ball of $\CC^n$.  See [KRA1, Ch.\ 11].  It is a relatively straightforward normal
families argument. We shall not repeat the details.
\endpf 
\smallskip \\

It is no longer the case in general (see our Proposition 3) that $\widehat{\M}$ equals 1 at one
point if and only if $\M$ equals 1 at all points---{\it unless} the model domain ${\cal B}$ has
transitive automorphism group.   See more on this point in what follows.

Now of course one of the great classical applications of Proposition 11, when ${\cal B}$ is the unit ball $B$,
is to prove the Bun Wong/Rosay theorem (our Theorem 2 above).   One might now ask whether a similar
sort of result could be proved with the new quotient invariant $\widehat{\M}$.   The answer is that
the {\it proof} requires that the model domain have transitive automorphism group (see the details
in [KRA1, Ch.\ 11]).

Thus we may only consider models ${\cal B}$ chosen from among the bounded
symmetric domains of Cartan (see [HEL]).   Let us concentrate here on the case
when ${\cal B}$ is the unit polydisc.	 The following result is similar to one
proved in [KIM]:

\begin{theorem} \sl
Let $\Omega \ss \CC^2$ be a smoothly bounded, convex domain.  Let $P \in \partial \Omega$
and assume that $\partial \Omega$ in a neighborhood $U$ of $P$ coincides with a real
hyperplane in $\CC^n$.  In suitable local coordinates we may say that
$$
\partial \Omega \cap U = \{z \in U: \Re z_1 = 0\} \, .	    \eqno (\dagger)
$$
If $P$ is a boundary orbit accumulation point for $\Omega$ then $\Omega$ is biholomorphic
to the bidisc.
\end{theorem}
{\bf Sketch of Proof:}  The key fact in the proof of this result when $P$ is a strongly pseudoconvex
point (our Theorem 3) is that the geometry localizes at $P$.   This means that if $X \in \Omega$
and $\varphi_j$ are automorphisms of $\Omega$ such that $\varphi_j(X) \ra P$ then $\varphi_j$
converges uniformly on any compact set $K$ to $P$.

Such is not the case in our present situation.   But the automorphisms $\varphi_j$ and the point
$X$ still exist (by a classical lemma of H. Cartan [NAR]).  As indicated in line $(\dagger)$, assume
that the real normal direction at $P$ is the $\Re z_1$ direction.  If $K \ss \Omega$ is any compact set 
then we may compose $\varphi_j$ for $j$ large with a dilation in the tangential directions
$z_2, z_3, \dots, z_n$ to localize the geometry near $P$, just as in the classical case.
The rest of the proof goes through as in the classical case described in [KRA1].  Instead
of localizing to an image of the ball, one localizes to a bidisc.
\endpf
\smallskip \\

\begin{remark} \rm
In [KIM], K.-T. Kim uses the method of {\it scaling} to obtain his result.  This is a powerful
technique that has wide applicability in this subject (see [GKK], for instance).  The argument
that we sketch here is similar in spirit to scaling.
\end{remark}

Perhaps another point worth considering is stability results for the quotient invariant $\M$ (i.e., the
original invariant modeled on the unit ball $B$).  We have the following result:

\begin{theorem} \sl
Let $\Omega, \Omega_j \ss \CC^n$ be bounded domains with $C^2$ boundary and suppose
that $\Omega_j \ra \Omega$ in the $C^2$ topology on domains (see [GRK3], [GRK4] for this
concept).  Then
$$
\M_{\Omega_j} \ra \M_\Omega
$$
uniformly on compact subsets of $\Omega$ as $j \ra \infty$.
\end{theorem}
{\bf Proof:}  Simply use the Carath\'{e}odory and Kobayashi stability
 results established in [GRK5].
 \endpf
 \smallskip \\

\section{Extremal Discs and Chains for the Kobayashi Metric}

In the remarkable paper [LEM], L. Lempert shows that, on a convex domain
$\Omega \ss \CC^n$, the integrated Kobayashi distance on $\Omega$ may be
calculated using a Kobayashi chain of length one disc (see [KOB], [KRA1]
for the concept of Kobayashi chain). This is done as a prelude to
developing his profound theory of extremal discs on strongly convex
domains.

Lempert comments that such a result is not true for general pseudoconvex domains, and he
provides the following example:

\begin{example} \rm
Let 
$$
\Omega_\epsilon = \{(z,w) \in \CC^2: |z| < 2, |w| < 2, |zw| < \epsilon\} \, .
$$
Let $P = (1,0) \in \Omega_\epsilon$ and $Q = (0,1) \in \Omega_\epsilon$.  Then the Kobayashi one-disc
distance of $P$ to $Q$ tends to infinity as $\epsilon \ra 0^+$.  Just to be perfectly clear,
we note that the one-disc Kobayashi distance of two points $P$ and $Q$ in a domain $\Omega$
is defined to be
$$
{\bf d}(P, Q) = \inf \{\rho(\varphi(a), \varphi(b)):  \varphi: D \ra \Omega, \varphi \ \hbox{holomorphic},
\varphi(a) = P, \varphi(b) = Q\} \, ,
$$
where $\rho$ is the classical Poincar\'{e} metric on the disc $D$.

Lempert's reasoning in this example (private communication) is as follows:  Suppose not.  Then
there are mappings $\varphi_\epsilon: D \ra \Omega_\epsilon$ with $\varphi_\epsilon(a_\epsilon) = P$
and $\varphi_\epsilon(b_\epsilon) = Q$ and $\rho(a_\epsilon, b_\epsilon)$ bounded above
as $\epsilon \ra 0^+$.  Thus we have that $a_\epsilon$, $b_\epsilon$ remain in a compact
subset $K$ of $D$.  Passing to a normal limit (with Montel's theorem), we find
a holomorphic function $\varphi_0: D \ra \{(z,w): |z| \leq 2, |w| \leq 2, |z \cdot w| = 1\}$
and points $a_0, b_0 \in K$ such that $\varphi_0(a_0) = P$, $\varphi_0(b_0) = Q$.  Of course
this is impossible, since it must be that either the image of $\varphi_0$
lies in $\{(z,w): z = 0\}$ or in $\{(z,w): w = 0\}$.  

It is useful, and instructive, to have a more constructive means of seeing that
this example works.  We thank John E. McCarthy for the following argument.

Take 
$$ 
\varphi = (f_1, f_2): D \ra \Omega_\epsilon 
$$
holomorphic.  We assume that
\begin{itemize}
\item $\varphi(0) = (1,0)$;
\item $\varphi(r) = (0,1)$.
\end{itemize}
We shall show constructively that, as $\epsilon \ra 0^+$, it must follow
that $r \ra 1^-$.   This is equivalent to what is claimed for the
domains $\Omega_\epsilon$.

Now use the inner-outer factorization for holomorphic functions on the disc (see, for
example [HOF]) to write $f_1 = F_1 \cdot I_1$ and $f_2 = F_2 \cdot I_2$.  Here
each $F_j$ is outer and each $I_j$ is inner.  Since $|f_1 \cdot f_2| < \epsilon$, 
we may be sure that 
$$
|F_1 \cdot F_2| < \epsilon \, .	  \eqno (*)
$$
Now certainly
$$
|F_1(0)| \geq |f_1(0)| = 1
$$
and hence
$$
|F_2(0)| < \epsilon \, .
$$
Certainly $\log |F_1| + \log |F_2|$ is harmonic, and by line $(*)$, is is majorized
by $\log \epsilon$.

Let $h$ denote the harmonic function $\log |F_2|$.  We can be sure that
\begin{enumerate}
\item[{\bf (1)}]  $h \leq \log 2$;
\item[{\bf (2)}]  $h(0) \leq \log \epsilon$;
\item[{\bf (3)}]  $h(r) \geq 0$.
\end{enumerate}

Let $h^+$ be the positive part of $h$ and $h^-$ the negative part.  Of
course $h^+ \geq 0$ and $h^- \geq 0$.  Then the mean-value property for 
harmonic functions tells us that
$$
\frac{1}{2\pi} \int_0^{2\pi} h^+ (e^{i\theta}) \, d\theta - \avgint h^-(e^{i\theta}) \, d\theta = h(0) \leq \log \epsilon
$$
hence
$$
\avgint h^-(e^{i\theta}) \, d\theta \geq \avgint h^+ (e^{i\theta}) \, d\theta + \log \frac{1}{\epsilon} \, .    \eqno (**)
$$

Let $P_r(e^{i\theta})$ denote the Poisson kernel for the unit disc $D$.  Then
$$
h(r) = \avgint h(e^{i\theta}) P_r(e^{i\theta}) \, d\theta \, .
$$
But Harnack's inequalities tell us that
$$
\frac{1-r}{1+r} \leq P_r(e^{i\theta}) \leq \frac{1 + r}{1 - r} \, .
$$
As a result, using {\bf (3)} above, 
$$
0 \leq h(r) \leq \frac{1 + r}{1 - r} \cdot \avgint h^+(e^{i\theta}) \, d\theta - \frac{1 - r}{1 + r} \cdot \avgint h^-(e^{i\theta}) \, d\theta \, . \eqno (*{*}*)
$$
We conclude that
$$
0 \leq h(r) \leq \frac{1 + r}{1 - r} \cdot \log 2 - \frac{1 - r}{1 + r} \cdot \avgint h^-(e^{i\theta} \, d\theta \, .
$$

Therefore
\begin{eqnarray*}
\avgint h^-(e^{i\theta}) \, d\theta & \leq & \left ( \frac{1 + r}{1 - r} \right )^2 \avgint h^+(e^{i\theta}) \, d\theta \\
                                    & \leq & \left ( \frac{1 + r}{1 - r} \right )^2 \left [ \avgint h^-(e^{i\theta}) \, d\theta + \log \epsilon \right ] \, ,
\end{eqnarray*}
where we have use $(**)$ in the last inequality.
Now certainly
\begin{eqnarray*}
            \log \frac{1}{\epsilon} & \leq & |h(0)| \\
                                    & \leq & \avgint h^+(e^{i\theta}) \, d\theta \\
                                    & \leq & \avgint h^-(e^{i\theta}) \, d\theta + \log \epsilon  \\
                                    & \leq & \avgint h^-(e^{i\theta}) \, d\theta   \\
      	                            & \leq & \left ( \frac{1 + r}{1 - r} \right )^2 \avgint h^+(e^{i\theta}) \, d\theta \\
                                    & \leq & \left ( \frac{1 + r}{1 - r} \right )^2 \cdot \log 2 \, .
\end{eqnarray*}
As $\epsilon \ra 0^+$, this last inequality can only be true if $r \ra 1^-$.  That is what we wished to prove.
\end{example}
\vspace*{.2in}

There has been some interest, since Lempert's paper, in developing an analogous theory on strongly
pseudoconvex domains.  N. Sibony [SIB] has shown that certain aspects of such a program
are impossible.

It is natural to reason as follows:
\begin{itemize}
\item Near the boundary of a strongly pseudoconvex domain, the domain is well approximated by
the biholomorphic image of $B$, the unit ball.  It is easy to verify directly (or by invoking
Lempert) that Kobayashi distance on the ball can be realized with Kobayashi chains of length 1.
\item In the interior of the domain---away from the boundary---things should be trivial.  After
all, if $\Omega$ is strongly pseudoconvex and $P \in \Omega$ is in the interior---away from the
boundary---then the infinitesimal Kobayashi metric $F_K^\Omega(P, \xi)$ for one Euclidean unit vector $\xi$ ought
to be roughly the same as the infinitesimal Kobayashi metric $F_K^\Omega(P, \xi')$ for any other Euclidean
unit vector $\xi'$.  Also the Kobayashi metric on a compact subset $K$ of $\Omega$ is comparable to the
Euclidean metric.  So one should be able to check directly that chains in the interior behave like chains
for the Euclidean metric.
\end{itemize}

Unfortunately the expectation enunciated in the second bulleted item above is not true.

\begin{example} \rm
Let $N > 0$ be a large positive integer and set
$$
B_N = \{(z_1, z_2) \in \CC^2: |z_1|^2 + |z_2/N|^2 < 1\} \, .
$$
Of course $B_N$ is biholomorphic to the unit ball $B$ via the biholomorphism
\begin{eqnarray*}
\Psi: B & \longrightarrow & B_N \\
    (z_1, z_2) & \longmapsto &  (z_1, Nz_2) \, .
\end{eqnarray*}

And one calculates readily, using the mapping $\Psi$, that
$$
F_K^{B_N} ((0,0), (1,0)) = 1
$$
while
$$
F_K^{B_N} ((0,0), (0,1)) = N \, .
$$
So the two different infinitesimal Kobayashi metric measurements at the base point ${\bf 0} = (0,0)$---in
two different Euclidean unit directions---are very different.
\end{example}

Interestingly, the following contrasting result is true for the Carath\'{e}odory metric:

\begin{proposition} \sl
Let $\Omega$ be a fixed, bounded domain in $\CC^n$.   Let $K \ss \Omega$ be a fixed compact
subset.  There is a positive constant $C_0$ so that, if $P \in K$ and $\xi_1, \xi_2$ are
Euclidean unit vectors then
$$
  \| F_C^\Omega(P, \xi_1) - F_C^\Omega(P, \xi_2) \|  \leq  C_0 \, .
$$
\end{proposition}
{\bf Proof:}  Let $r > 0$ be a small number.   Let $\gamma$ be a $C_c^\infty$ function that satisfies:
\begin{enumerate}
\item[{\bf (a)}]  $\gamma$ is radial.
\item[{\bf (b)}]  $\gamma$ is supported in the Euclidean ball with center at $P$ and radius $r$.
\item[{\bf (c)}]  $\gamma$ is identically equal to 1 on the Euclidean ball with center at $P$ and radiuis $r/2$.
\end{enumerate}						       

Now let $\mu$ be a unitary rotation of $\CC^n$ that takes $\xi_1$ to $\xi_2$.   Fix a point $P \in K$ and
vectors $\xi_1, \xi_2$ as in the statement of the proposition.  Let $\psi$ be an element
of $(\Omega, D)$ with $\psi(P) = 0$ and $\psi'(P)$ a positive, real multiple of $\xi$---say that $\psi'(P) = \kappa \xi$.
Set 
$$
\widetilde{\psi}(z) = \gamma(z) \cdot [\psi \circ \mu^{-1}(z)] + [1 - \gamma(z)] \cdot \psi(z) + h(z) \, .  \eqno (*)
$$
Of course $\widetilde{\psi}$ will not be {\it a priori} holomorphic---because we have constructed the
function using cutoff functions---but we hope to use the $\overline{\partial}$ problem to select $h$
so that $\widetilde{\psi}$ {\it will} be holomorphic.

Applying the $\overline{\partial}$ operator to both sides of equation $(*)$, we find
that
$$
\overline{\partial} h = - \dbar \gamma \cdot [\psi \circ \mu^{-1}] + \dbar \gamma \cdot \psi \, .
$$
Now it is essential to notice these properties:
\begin{itemize}
\item $|\dbar \gamma|$ is of size $\approx 1/r$;
\item $\psi(P) = 0$, so that, on the support of $\dbar \gamma$, $|\psi|$ of size $r$;
\item $\dbar h$ is supported on the ball with center $P$ and radius $r$;
\item $\dbar h$ is $\dbar$-closed.
\end{itemize}
We see therefore that $\dbar h$ is of size ${\cal O}(1)$ (in Landau's notation) and supported
in a Euclidean ball of radius $r$.  Hence $h$ has $L^3$ norm on any one-dimensional complex
slice of space not exceeding  $C \cdot [r^2 \cdot 1]^{1/3} = C \cdot r^{2/3}$.

Now we may solve the equation $\dbar u = h$ using the solution 
$$
h(z) = - \frac{1}{\pi} \int \! \! \int \frac{\tau_j(z_1, \dots, z_{j-1}, \xi, z_{j+1}, \dots, z_n)}{\xi - z_j} \, dA(\xi) \, .
$$
Here $\dbar h = \tau_1 d\overline{z}_1 + \tau_2 d\overline{z}_2$.  [See [KRA1, p.\ 16] for a discussion of this idea.]   Then we see that
$$
\|u\|_{\rm sup} \leq \left \| \tau_j(z_1, \dots, z_{j-1}, \, \cdot \, , z_{j+1}, \dots, z_n) \right \|_{L^3} \cdot \left \| \frac{1}{ \, \cdot \, - z} \right \|_{L^{3/2}}
	       \leq r^{2/3} \cdot r^{1/2} = r^{7/6} \, .
$$
In summary, $h$ is small in uniform norm if $r$ is small, and we may choose $r$ in advance to be as small as we please.

Now what is more essential for our purposes is that we may likewise estimate the size of $\|\nabla h\|_{\rm sup}$.  For we may 
write
$$
h(z) = - \frac{1}{\pi} \int \! \! \int \frac{\tau_j(z_1, \dots, z_{j-1}, \xi - z_j, z_{j+1}, \dots, z_n)}{\xi} \, dA(\xi) 
$$
and hence
$$
\nabla h(z) = - \frac{1}{\pi} \int \! \! \int \frac{\nabla_z \tau_j(z_1, \dots, z_{j-1}, \xi - z_j, z_{j+1}, \dots, z_n)}{\xi} \, dA(\xi)  \eqno (\ddagger)
$$
But now it is essential to notice that 
\begin{itemize}
\item $|\nabla \dbar \gamma|$ is of size $r^{-2}$;
\item $\nabla \psi$ is of size ${\cal O}(1)$.
\end{itemize}
It follows then that $\nabla \dbar h$ is of size $r^{-1}$ and is still supported on a Euclidean ball of radius $r$.
Thus we may estimate $(\ddagger)$ again using H\"{o}lder's inequality.  The result is that
$\|\nabla h\|_{\rm sup} \leq C \cdot r^{1/6}$.

We conclude that the corrected candidate function $\widetilde{\psi}$ is, near $P$ uniformly closed to being just a rotation
of $\psi$.   We also see that
$$
\widetilde{\psi}'(P) = \psi'(P) \circ \mu + h'(P) \, .
$$
Thus $\widetilde{\psi}'(P)$ is as close as we like to equalling $\xi'$.  
Now taking a normal limit (again using Montel's theorem) as $r \ra 0^+$ yields a function $\psi_0: \Omega \ra B$
with $\psi_0(P) = 0$ and $\psi'_0(P) = \kappa \cdot \xi'$.   So we find a candidate for the Carath\'{e}odory metric
at $P$ in the direction $\xi'$ that is comparable to the original candidate $\psi$ in the direction $\xi$.
\endpf 
\smallskip \\

We would like to explore here the nature of Kobayashi chains on a strongly pseudoconvex domain.
In principle, the Kobayashi chains on a given domain $\Omega$ could have any number
of discs.  We shall prove, however, that on a strongly pseudoconvex domain there
is an {\it a priori} upper bound for the length of chains.  This result may be thought
of as a prelude to the development of a Lempert-type theory on strongly pseudoconvex domains.

\begin{proposition} \sl Let $\Omega \ss \CC^n$ be a strongly pseudoconvex domain with $C^2$ boundary.
Let $f: D \ra \Omega$ and $g: D \ra \Omega$ be holomorphic mappings of the disc
into $\Omega$.  We assume that $\sup_{\zeta \in D} |\varphi_1(\zeta) - \varphi_2(\zeta)| < \delta$
for some small $\delta > 0$.   Further, following Lempert's notation [LEM, pp.\ 430--431], we
let $\zeta, \omega, \omega' \in D$ satisfy
$$
f(\zeta) = z \, , \quad f(\omega) = g(\omega') = w \, , \quad g(\sigma) = s \, .
$$
Then there is a holomorphic mapping 
$$
h: D \ra \Omega
$$
with $h(\zeta) = z$, $h(\sigma) = s$.  It follows then that, in the calculation of the 
Kobayashi metric using chains, we may replace the two discs $f$, $g$ with the single
disc $h$.
\end{proposition}
{\bf Proof:}  By the Forn\ae ss imbedding theorem, there is a strongly
convex domain $\Omega'$ with $C^2$ boundary, $\Omega' \ss \CC^N$ with
$N > \, > n$ in general, and a proper holomorphic imbedding
$$
\Phi: \overline{\Omega} \ra \overline{\Omega'} \, .
$$
We refer the reader to [FOR] for the details of the domain
and the mapping.  Let $\widehat{\Omega} \ss \Omega'$ be the 
image of $\Omega$ under the mapping $\Phi$.  According to the Docquier-Grauert
theorem ([ROS], [DOG]), there is a neighborhood $U$ of $\widehat{\Omega}$ and
a holomorphic retraction $\pi: U \ra \widehat{\Omega}$.

Of course $\Phi(f(D))$ and $\Phi(g(D))$ both lie in $\widehat{\Omega}$. We
may apply Lempert's Theorem 1 to obtain a convex combination $\lambda(\zeta)$ of
$\Phi(f(D))$ and $\Phi(g(D))$.  Now we may not conclude that the image of $\lambda$
lies in $\widehat{\Omega}$.  But it certainly lies in the strongly convex domain $\Omega'$.
And, if $\delta$ is sufficiently small, then we know that the image of $\lambda$ lies in $U$.
Thus we may consider the analytic disc $\widehat{\lambda} \equiv \pi \circ \lambda$, whose image {\it does}
lies in $\Omega'$.   Now $\Phi^{-1}$ makes sense on $\Omega'$, so we may define
$$
h(\zeta) = \Phi^{-1} \circ \widehat{\lambda} \, .
$$
Tracing through the logic shows that this $h$ is the one that we seek.
\endpf
\smallskip \\

\begin{theorem} \sl Let $\Omega \ss \CC^n$ be a strongly pseudoconvex
domain with $C^2$ boundary. Then there is an $\epsilon > 0$ and an {\it a
priori} constant $K = K(\Omega)$ so that if $P, Q \in \Omega$ then there
is a Kobayashi chain with elements $\varphi_1$, \dots, $\varphi_k$ so that
the integrated Kobayashi distance of $P$ to $Q$ is within $\epsilon$ of
the length given by the Kobayashi chain.
\end{theorem}
{\bf Proof:}  Since $\Omega$ is a bounded domain, it is contained in a large Euclidean
 ball.  By elementary comparisons, (see [KRA1]), we know that the Kobayashi metric
 in $\Omega$ is not less than the Kobayashi metric in the ball.   In particular, we get
 an {\it a priori} upper bound on derivatives of extremal discs for the Kobayashi metric
 in $\Omega$.  As a result, there is an $\eta > 0$ and a finite net of points ${\cal P} \ss \Omega$
 so that 
\begin{enumerate}
\item[{\bf (i)}]  Every point of $\Omega$ is Euclidean distance not more than $\eta$ from
some point of ${\cal P}$;
\item[{\bf (ii)}]  There is an {\it a priori} integer $M > 0$ so that if $\psi: D \ra \Omega$ is
a Kobayashi extremal disc then there is a collection of elements ${\cal Q}_{\psi}$ of
at most $M$ points in ${\cal P}$ so that every point in the image $\psi(D)$ is Euclidean distance
at most $\eta$ from some point of ${\cal Q}_\psi$.   More importantly, there is a finite
net of points ${\cal K}_\psi$ in the disc $D$---of cardinality at most $M$---so that every
element of ${\cal Q}_\psi$ is the approximate image (within distance $\eta$) under $\psi$ of some element of ${\cal K}_\psi$
(in fact one can conveniently take ${\cal K}_\psi$ to be a net in the disc $D$ that has
unit distance $\eta'$, for some small $\eta' > 0$, in the Poincar\'{e} metric).  Thus
we associate to $\psi$ the set ${\cal Q}_\psi^{{\cal K}_\psi}$.
\end{enumerate}

Of course there are only finitely many possible sets ${\cal K}_\psi {\cal Q}_\psi$ (indeed $2^M$ is an
upper bound on the cardinality of $\{ {\cal Q}_\psi\}$, and there is a similar upper bound $2^{M'}$
for the $\{{\cal K}_\psi\}$).  If ${\cal T}$ is a Kobayashi
chain in $\Omega$ with more than ${2^M}^{2^{M'}}$ discs, then two of those discs will
share the same ${\cal K}_\psi$ and ${\cal Q}_\psi$.  As a result, if $\eta$ and $\eta'$ are fixed small enough (depending
on $\delta$ in the last proposition), then the two corresponding extremal discs in the chain will be close
enough that the last proposition applies.  And those two discs may be replaced by a single disc.

This shows that our {\it a priori} constant $K$ exists and does not exceed ${2^M}^{2^{M'}}$.
\endpf 
\smallskip \\

\section{Concluding Remarks}

In the past forty years or more, the Carath\'{e}odory and Kobayashi metric
constructions have proved to be powerful tools in both geometry and function
theory.  Their role in the study of automorphism group is more recent, but
is equally significant.  We trust that the contributions of this paper will point
in some new directions in the subject.   What lies in the future can only
be a topic for omphaloskepsis.

\newpage

\noindent {\Large \sc References}
\medskip  \\

\begin{enumerate}

\item[{\bf [BEL]}] S. Bell, Biholomorphic mappings and the $\overline
\partial $-problem, {\it Ann.\ of Math.} 114(1981), 103--113.

\item[{\bf [CAT]}] D. Catlin, Estimates of invariant metrics on
pseudoconvex domains of dimension two, {\em Math. Z.} 200(1989), 429-466.

\item[{\bf [EIS]}] D. Eisenman, {\it Intrinsic Measures on Complex
Manifolds and Holomorphic Mappings}, Memoir of the American Mathematical
Society, Providence, RI, 1970.

\item[{\bf [FOR]}] J. E. Forn\ae ss, Strictly pseudoconvex domains in
convex domains, {\em Am. J. Math.} 98(1976), 529-569.

\item[{\bf [GKK]}] R. E. Greene, K.-T. Kim, and S. G. Krantz, {\it The
Geometry of Complex Domains}, Birkh\"{a}user Publishing, Boston, MA, 2010,
to appear.

\item[{\bf [GRK1]}] R. E. Greene and S. G. Krantz, Biholomorphic self-maps
of domains, {\it Complex Analysis II} (C. Berenstein, ed.), Springer
Lecture Notes, vol. 1276, 1987, 136-207.

\item[{\bf [GRK2]}] R. E. Greene and S. G. Krantz, Techniques for Studying
the Automorphism Groups of Weakly Pseudoconvex Domains, Proceedings of the
Special Year at the Mittag-Leffler Institute (J. E. Forn\ae ss and C. O.
Kiselman, eds.) {\it Annals of Math. Studies,} Princeton Univ. Press,
Princeton, 1992.

\item[{\bf [GRK3]}] R. E. Greene and S. G. Krantz, Stability properties of
the Bergman kernel and curvature properties of bounded domains, {\it
Recent Developments in Several Complex Variables} (J. E. Forn\ae ss, ed.),
Princeton University Press (1979), 179-198.

\item[{\bf [GRK4]}] R. E. Greene and S. G. Krantz, Deformations of complex
structure, estimates for the $\dbar$-equation, and stability of the
Bergman kernel, {\it Advances in Math.} 43(1982), 1-86.

\item[{\bf [GRK5]}] R. E. Greene and S. G. Krantz, Stability of the
Carath\'{e}odory and Kobayashi metrics and applications to biholomorphic
mappings, {\it Proceedings of Symposia in Pure Mathematics,} vol. 41
(1984), 77-93.						      

\item[{\bf [HEL]}]  S. Helgason, {\it Differential Geometry and Symmetric Spaces},
Academic Press, New York, 1962.

\item[{\bf [HOF]}] K. Hoffman, {\it Banach Spaces of Holomorphic
Functions}, Prentice-Hall, Englewood Cliffs, 1962.
	
\item[{\bf [ISK]}] A. Isaev and S. G. Krantz, Domains with
non-compact automorphism group: A Survey, {\it Advances in
Math.} 146(1999), 1--38.

\item[{\bf [KIM]}] K.-T. Kim, Domains in $C^n$ with a piecewise
Levi flat boundary which possess a noncompact automorphism
group, {\it Math.\ Ann.} 292(1992), 575--586.

\item[{\bf [KOB]}]  S. Kobayashi, {\it Hyperbolic Manifolds and Holomorphic
Mappings}, Dekker, New York, 1970.

\item[{\bf [KRA1]}]  S. G. Krantz, {\it Function Theory of Several Complex
Variables}, $2^{\rm nd}$ ed., American Mathematical Society, Providence,
RI, 2001.

\item[{\bf [KRA2]}] S. G. Krantz, The Carath\'{e}odory and Kobayashi
metrics and applications in complex analysis, {\it American Mathematical
Monthly} 115(2008), 304--329. 

\item[{\bf [KRA3]}] S. G. Krantz, The boundary behavior of the Kobayashi
metric, {\it Rocky Mountain Journal of Mathematics} 22(1992), 227--233.

\item[{\bf [LEM]}] L. Lempert, La metrique \mbox{K}obayashi et las
representation des domains sur la boule, {\em Bull. Soc. Math. France}
109(1981), 427-474.

\item[{\bf [MCN]}] J. McNeal, Convex domains of finite type,
{\it J. Funct.\ Anal.} 108(1992), 361--373.

\item[{\bf [NAR]}]  R. Narasimhan, {\it Several Complex Variables}, University
of Chicago Press, Chicago, 1971.
				       
\item[{\bf [RUD]}] W. Rudin, {\it Function Theory in the Unit Ball of
$\CC^n$}, Grundlehren der Mathematischen Wissenschaften in
Einzeldarstellungen, Springer, Berlin, 1980.

\item[{\bf [SIB]}]  N. Sibony, unpublished notes.

\item[{\bf [WEB]}] S. M. Webster, On the mapping problem for algebraic real
hypersurfaces, {\it Invent.\ Math.} 43(1977), 53--68.

\item[{\bf [WON]}] B. Wong, Characterizations of the ball in $\CC^n$ by its
automorphism group, {\em Invent. Math.} 41(1977), 253-257.

\end{enumerate}
\vspace*{.4in}

\begin{quote}
Department of Mathematics \\
Washington University in St.\ Louis \\
St.\ Louis, Missouri 63130 \ \ USA  \\
{\tt sk@math.wustl.edu}  \\
\end{quote}

\end{document}